\documentclass[12pt,reqno]{amsart}
\usepackage[cp1251]{inputenc}
\usepackage[english]{babel}
\usepackage[dvips]{graphicx}
\usepackage{amsmath,amsopn,amssymb,amsthm}

\newtheorem{theorem}{Theorem}
\newtheorem{remark}{Remark}

\begin{document}

\title
{Generalized Popoviciu’s problem}
\author{Yu.G. Nikonorov, Yu.V. Nikonorova}

\begin{abstract}
This is an English translation of the following paper, published several years ago:
{\it
Nikonorov Yu.G., Nikonorova Yu.V.
Generalized Popoviciu’s problem (Russian),
Tr. Rubtsovsk. Ind. Inst., 7, 229--232 (2000), Zbl. 0958.51021}.
All inserted footnotes provide additional information related to the mentioned problem.

\vspace{2mm} \noindent 2010 Mathematical Subject Classification:
51M04, 51M16, 52A10, 52A40.

\vspace{2mm} \noindent Key words and phrases: Popoviciu’s problem, generalized Popoviciu’s problem, convex quadrilaterals, extremal problems.
\end{abstract}

\maketitle

Let $ABCD$ be any convex quadrilateral in the Euclidean plane.
Let us consider the points
$A_1$, $B_1$, $C_1$, and $D_1$ of the segments
 $AB$, $BC$, $CD$, and $DA$ respectively, such that
$$
\frac{|A A_1|}{|A_1B|}=
\frac{|B B_1|}{|B_1C|}=
\frac{|C C_1|}{|C_1D|}=
\frac{|D D_1|}{|D_1A|}=k
$$
for some fixed $k>0$.
The straight lines
$AB_1$, $BC_1$, $CD_1$, $DA_1$ form a quadrilateral $KLMN$
($K$, $L$, $M$, $N$ are the intersection points for the first and the fourth,
the first and the second, the second and the third, the third and the fourth straight lines respectively),
situated inside $ABCD$.

\begin{center}
\begin{figure}[th]
\centering\scalebox{1}[1]{\includegraphics[angle=0,totalheight=3in]
{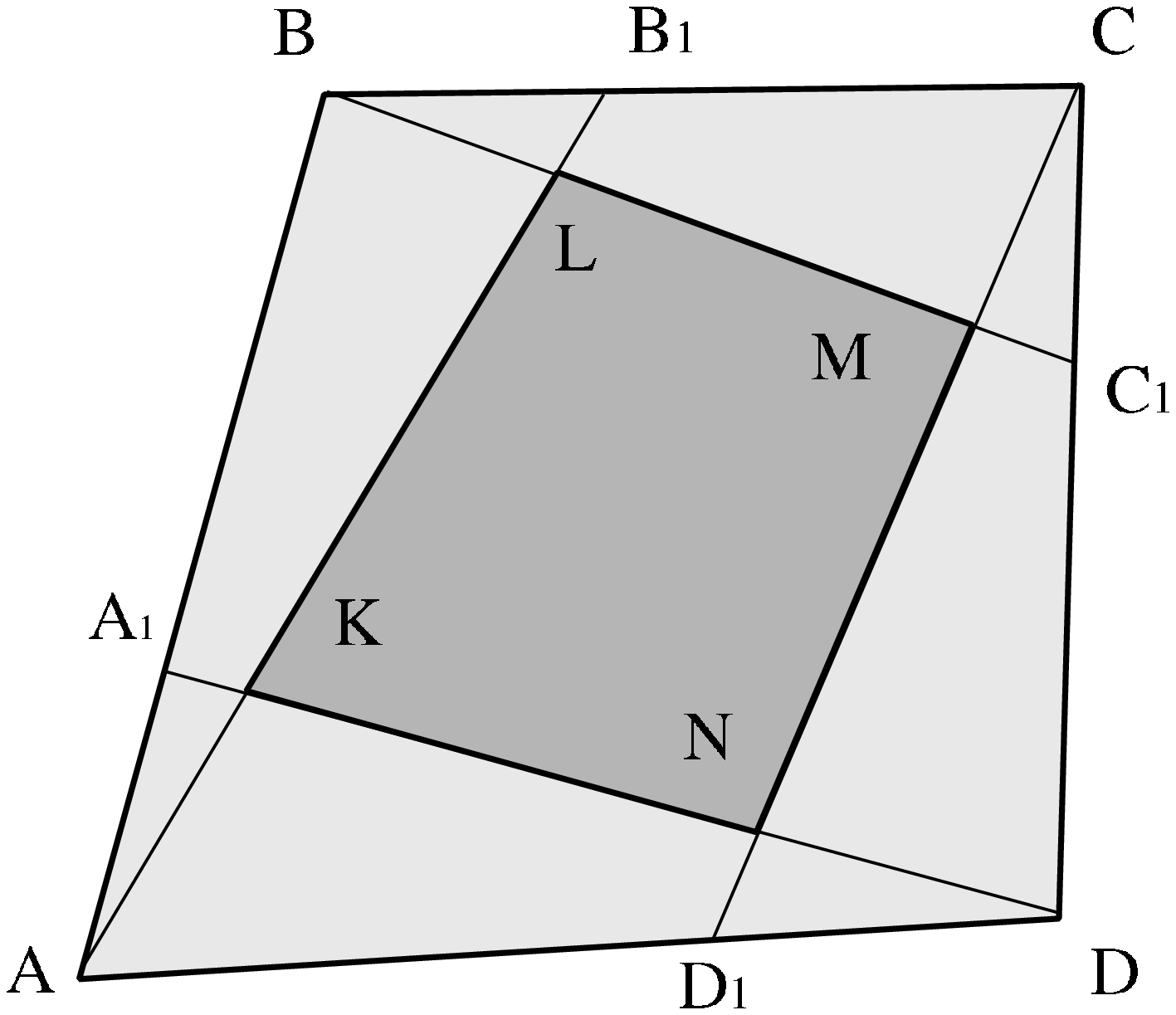}}
\caption{}\label{popov}
\end{figure}
\end{center}

Let us denote by $S$ and $s$ the areas of the quadrilaterals
$ABCD$ and $KLMN$ respectively. Popoviciu’s problem
\footnote{\,T.~Popoviciu, Problem 5897, Gazeta Matematic\u{a}, 49 (1943), P.~322; see also \\ \url{http://www.gazetamatematica.net/?q=node/2135}},
(see e.~g., \cite[P.~39]{Shkl}
\footnote{\,See also pp.~132--133 in the following book: {\it Bottema~O., {\DJ}or{\dj}evi\'c~R.\v{Z}., {Jani\'c} R.R., {Mitrinovi\'c}~D.S., {Vasi\'c}~P.M.
Geometric Inequalities, Wolters-Noordhoff Publishing, Groningen, 1969, 151~p.}})
consists in proving the inequality
$$
\frac{1}{6}\,S \leq s \leq \frac{1}{5}\,S
$$
for $k=1$ and in the study of cases of equalities.
The solution of original Popoviciu’s problem was obtained by
Yu.G.~Nikonorov in \cite{Nikk}.
\footnote{\,Later the same result was obtained also in the following papers:
{\it Ash~J.M, Ash~M.A, Ash~P.F.
Constructing a quadrilateral inside another one, The Mathematical Gazette, 93(528), 522--532 (2009);
Mabry~R.
Crosscut convex quadrilaterals,
Math. Mag., 84(1), 16--25 (2011), Zbl~1227.51015}.
The following papers are devoted to a study of the corresponding analogue for pentagons on the Euclidean plane:
{\it Dudkin F.A.
Popoviciu’s problem for a convex pentagon (Russian), Tr. Rubtsovsk. Ind. Inst., 12, 31--38 (2003), Zbl 1036.52002;
Kizbikenov~K.O. Popoviciu’s problem for a convex pentagon (Russian),
Vestnik Altaiskogo Gosudarstvennogo Pedagogiceskogo Universiteta, 20, 16--23 (2014), \url{https://elibrary.ru/item.asp?id=21991088}}.}
In this paper we obtain a more general result that valid for all $k>0$.

\begin{theorem}
\label{cr2.1.thm.1}
For an arbitrary convex quadrilateral in the Euclidean plane and for any $k>0$,
the inequality
$$
\frac{1}{(k+1)(k^2+k+1)}\,S \leq s \leq \frac{1}{2k^2+2k+1}\,S\,
$$
holds. Moreover, if $S>0$ then the equality $(k+1)(k^2+k+1)s=S$ is fulfilled exactly for
quadrilaterals with two coinciding vertices, whereas any quadrilateral with the property $(2k^2+2k+1)s=S$
is situated in a continuous family of quadrilaterals with the same property,
which contains a parallelogram.\footnote{\,The same result was obtained later in the following paper:
{\it Ash~J.M, Ash~M.A, Ash~P.F.
Constructing a quadrilateral inside another one, The Mathematical Gazette, 93(528), 522--532 (2009)}, see also arXiv:0704.2716.}
\end{theorem}

We will give a proof that does not depend on the use of computer technology. Nevertheless, it is difficult to obtain manually
the factorization of polynomials, that are used in the proof.
All the relevant calculations can easily be done with the help  of symbolic computation systems.
\footnote{\,See e.~g., Application 1 in the following book:
{\it Nikonorov Yu.G., Nikonorova, Yu.V.
Applications of MAPLE to the solution of geometric problems. Textbook. 2-nd ed. (Russian).
Rubtsovsk: Izdatel’stvo Altajskogo Gosudarstvennogo Universiteta (ISBN 5-7904-0290-9/pbk). 78 p. (2005), Zbl 1095.65017}.
}

\begin{proof}
We may assume that the vertices $A$, $B$, and $D$ are pairwise distinct points.
Using a suitable affine transformation
(the ratio of the areas does not change in this case), one can reduce the problem to the case,
when $\angle BAD$ of the quadrilateral $ABCD$ is right,
and the sides $AB$ and $AD$ have unit length.
We introduce a Cartesian coordinate system in the plane,
taking the point $A$ as the origin and the rays $AD$
and $AB$ as the coordinate rays. In this coordinate system, the points $A$, $B$, $D$, $C$
have coordinates $(0,0)$, $(0,1)$, $(1,0)$, $(a,b)$ respectively,
where $a \geq 0$, $b \geq 0$, $a+b \geq 1$. Let us denote by $\Omega$
the set
$$
\{(a,b)\in \mathbb{R}^2 \,\,|\,\,a \geq 0, \, b \geq 0, \, a+b \geq 1 \}\, .
$$
Now, let us calculate the values $s$ and $S$
using standard analytical
geometry tools.
It is easy to get that $2S=a+b$.
It is a simple problem to calculate the coordinates of the following points:
\begin{eqnarray*}
A_1&=&\left(0,\frac{k}{k+1}\right),\\
B_1&=&\left(\frac{ak}{k+1}, \frac{1+kb}{k+1}\right),\\
C_1&=&\left(\frac{a+k}{k+1}, \frac{b}{k+1}\right),\\
D_1&=&\left(\frac{1}{k+1},0\right).
\end{eqnarray*}
It is easy also to find the equations of the straight lines $DA_1$, $CD_1$, $AB_1$, and $BC_1$:
\begin{eqnarray*}
kx+(k+1)y-k=0,
\\
(k+1)bx+(1-a(k+1))y-b =0,
\\
(kb+1)x-kay =0,
\\
(k+1-b)x+(a+k)y-(a+k) =0.
\end{eqnarray*}

Now, we should calculate the coordinates of the points $K$, $L$, $M$, and $N$, which are intersection points
of pairs of the corresponding straight lines. Omitting the standard calculations, we obtain
\renewcommand{\arraystretch}{1.4}
{\small
\begin{eqnarray*}
K=\left(\frac{ak^2}{ak^{2}+bk^{2}+bk+k+1},\,
 \frac{k(bk+1)}{ak^{2}+bk^{2}+bk+k+1} \right),
\\
L=\left(\frac{ak(a+k)}{ak^{2}+bk^{2}+ak+k+a},\,
\frac{(bk+1)(a+k)}{ak^{2}+bk^{2}+ak+k+a} \right),
\\
M=\left(\frac{ak^{2}+a^{2}k+ak+bk-k+a^{2}+ab-a}{ak^{2}+bk^{2}+2ak+bk-k+a+b-1},\,
\frac{b\,(k^{2}+ak+a+b-1)}{ak^{2}+bk^{2}+2ak+bk-k+a+b-1} \right),
\\
N= \left(
\frac{ak^{2}+ak+bk-k+b}{ak^{2}+bk^{2}+ak+2bk-k+b},\,
\frac{bk^{2}}{ak^{2}+bk^{2}+ak+2bk-k+b}\right).
\end{eqnarray*}
}
It is obvious that
$$
2s=2S_{\Delta ANM}+2S_{\Delta AML}-2S_{\Delta ANK}.
$$

Using determinants for calculating the areas of the triangles, we obtain

{\small
\begin{eqnarray*}
2S_{\Delta AML}=
\left( a+k\right)\cdot \frac{ak^{2}-bk^{2}-k+ak^{2}b+ak+a^{2}k+b^{2}k^{2}+
bk-a+ab+a^{2}}
{\left(-k+ak^{2}+2ak-1+a+b+bk^{2}+bk\right)
 \left( ak^{2}+ak+a+bk^{2}+k\right)},
\\
2S_{\Delta ANM}=
b\,\cdot \frac{ak^{2}b+k-b-2bk+ab-2ak-ak^{2}+bk^{2}+b^{2}k+a^{2}k+b^{2}+3bka+
a^{2}k^{2}}{\left( bk^{2}+2bk+b-k+ak^{2}+ak\right)
\left(-k+ak^{2}+2ak-1+a+b+bk^{2}+bk\right) },
\\
2S_{\Delta ANK}=
k\,\cdot \frac{b^{2}k^{2}+b^{2}k-bk^{2}+ak^{2}b+bk+b-k+ak^{2}+ak}
{\left(bk^{2}+bk+k+1+ak^{2}\right)
\left( bk^{2}+2bk+b-k+ak^{2}+ak\right) }.
\end{eqnarray*}
}
Therefore,
$$
\frac{s}{S}=\frac{P(a,b)}{Q(a,b)},
$$
where
\begin{eqnarray*}
P(a,b)=-2a^{2}k^{2}b+6ab^{4}k^{3}-4ak^{2}b+12a^{2}k^{3}b^{2}+9ab^{4}k^{5}
\\
+16a^{2}k^{4}b^{2}+19a^{2}k^{4}b^{3}+9a^{2}b^{3}k^{2}+8a^{3}k^{2}b^{2}+17a^{3}k^{3}b^{2}+2a^{2}k^{6}b
\\
+2a^{4}k^{6}b+8a^{3}k^{5}b-b^{2}a-3b^{2}k^{2}+a^{2}k+bk^{2}+ak^{2}-8ab^{3}k^{5}+9a^{4}k^{4}b
\\
+2b^{3}ak-4a^{2}k^{5}b+6a^{3}kb+3b^{4}k^{3}-5a^{2}kb+3ab^{4}k^{2}+14a^{3}k^{5}b^{2}
\\
+18a^{2}k^{5}b^{3}-4a^{2}k^{5}b^{2}+6a^{3}k^{6}b^{2}+15a^{2}k^{3}b^{3}+11ab^{4}k^{4}+ab^{2}k^{2}-9ak^{3}b^{2}
\\
+12a^{3}k^{4}b+4ab^{3}k^{4}+8a^{3}k^{2}b+a^{4}k^{3}+3a^{4}k^{2}-4a^{2}k^{6}b^{2}-3a^{3}k^{2}+8a^{4}k^{3}b
\\
+a^{3}k^{3}-2a^{2}k^{3}-5a^{3}k^{4}+6a^{2}k^{6}b^{3}+2b^{2}k^{3}-a^{2}b-2a^{3}k+2ab^{4}k^{6}
\\
+17a^{3}k^{4}b^{2}+5a^{4}k^{5}b+a^{4}k^{4}+a^{5}k^{5}+4a^{3}k^{3}b+12ab^{3}k^{3}+4ab^{3}k^{2}
\\
+7a^{2}k^{2}b^{2}-2b^{2}ka+7a^{2}b^{2}k+8ak^{4}b+8ak^{5}b^{2}-5b^{3}k^{3}+b^{3}k^{4}+4b^{3}k^{5}
\\
-2bk^{4}+2b^{2}k^{4}+6a^{2}k^{4}-4b^{2}k^{5}-2ak^{4}+4a^{2}k^{5}+a^{4}k+2a^{2}b^{2}+a^{3}b
\\
+2ak^{6}b^{2}-11ak^{4}b^{2}-3a^{2}k^{3}b-17a^{2}k^{4}b-b^{2}k+2a^{4}k^{6}+4a^{4}k^{5}-8a^{3}k^{5}
\\
+2a^{5}k^{4}+2b^{4}k^{6}-b^{4}k^{4}-a^{2}k^{2}+a^{5}k^{3}+a^{3}b^{2}k-2a^{3}k^{6}+2a^{2}b^{3}k
\\
-2b^{3}k^{6}+2a^{4}k^{2}b+b^{4}ak+b^{3}a+b^{5}k^{3}+2b^{4}k^{2}+b^{3}k+2b^{5}k^{4}+b^{5}k^{5},
\end{eqnarray*}

\begin{eqnarray*}
Q(a,b)=(a+b)\left( ak^{2}+ak+a+bk^{2}+k\right)
\\
\times \left(-k+ak^{2}+2ak-1+a+b+bk^{2}+bk\right)
\\
\times \left( bk^{2}+bk+k+1+ak^{2}\right)
\left( bk^{2}+2bk+b-k+ak^{2}+ak\right).
\end{eqnarray*}

It should be noted that
\begin{eqnarray*}
-k+ak^{2}+2ak-1+a+b+bk^{2}+bk
\\
=(a+b-1)+k(a+b-1)+ak+k^2(a+b)>0,
\\
bk^{2}+2bk+b-k+ak^{2}+ak=k(a+b-1)+bk+b+k^2(a+b)>0
\end{eqnarray*}
on the set $\Omega$. This implies that $Q(a,b)>0$ on the set $\Omega$.
\smallskip

Therefore, the proof of the theorem reduces to the proof of the following two inequalities for $(a,b)\in \Omega$:
\smallskip

1) $(k+1)(k^2+k+1)P(a,b)-Q(a,b) \geq 0$,

\smallskip

2) $Q(a,b)-(2k^2+2k+1)P(a,b) \geq 0$,
\smallskip

\noindent and the study of the equality cases.
\smallskip

Let us consider the first inequality.
Direct calculations show that

\begin{eqnarray*}
(k+1)(k^2+k+1)P(a,b)-Q(a,b)=k^{3} \left( a+b\right)\left( 1+2k+2k^{2}\right)
\\
\times \left(ak^{2}b+ak^{2}+2ak+a+b-1-k+b^{2}k-bk^{2}+b^{2}k^{2}\right)
\\
\times \left( a^{2}k^{2}+a^{2}k+ba-2ak+k+2bak-ak^{2}+ak^{2}b+bk^{2}\right).
\end{eqnarray*}

Note that
\begin{eqnarray*}
ak^{2}b+ak^{2}+2ak+a+b-1-k+b^{2}k-bk^{2}+b^{2}k^{2}
\\
=(a+b-1)(k^2b+1)+ak^2+(2a+b^2-1)k,
\\
a^{2}k^{2}+a^{2}k+ba-2ak+k+2bak-ak^{2}+ak^{2}b+bk^{2}
\\
=ak^2(a+b-1)+ab+bk^2+(a^2+2ab-2a+1)k.
\end{eqnarray*}

We note that  the inequality $a+b \geq 1$ is fulfilled on the set $\Omega$.
Thus, to prove the first part of the theorem it suffices
to make sure that the inequalities
$2a+b^2-1 \geq 0$
and
$a^2+2ab-2a+1 \geq 0$
are fulfilled on
the set $\Omega$.

Since the straight line $a+b=1$ is tangent to the convex curve
$2a+b^2-1=0$ (that is a parabola) at the point $(0,1)$, then it suffices to verify that
the inequality $2a+b^2-1 > 0$ is fulfilled at least at one point of the set $\Omega$.
It is easily to check it at the point $(1,1)$.
Notice, that
$2a+b^2-1= 0$  and
$(a+b-1)(k^2b+1)+ak^2+(2a+b^2-1)k=0$ simultaneously exactly for $(a,b)=(0,1)$.

Obviously,
the inequality $a^2+2ab-2a+1 =(a-1)^2+2ab \geq 0$ is fulfilled on the set
$\Omega$, and the equality is achieved only when
$(a,b)=(1,0)$. Note that for these values of variables,
the equality
$ak^2(a+b-1)+ab+bk^2+(a^2+2ab-2a+1)k=0$ holds.

Finally, we get that
$$
(k+1)(k^2+k+1)P(a,b)-Q(a,b) \geq 0
$$
on the set $\Omega$, and the equality is achieved only when
$(a,b)=(1,0)$ or $(a,b)=(0,1)$, that is, in the case when two vertices of
the quadrilateral are coincided.
\smallskip

The proof of the second part of the theorem is much simpler.
By direct calculations we obtain that
\begin{eqnarray*}
Q(a,b)-(2k^2+2k+1)P(a,b)
\\
=k^{4}\left( a+b\right)
\left( bk+1-a-ak\right) ^{2}
\left(b+bk-1+ak-2k\right) ^{2}.
\end{eqnarray*}

We note that the equality holds on the straight lines
$$
bk+1-a-ak=0,\qquad
b+bk-1+ak-2k=0,
$$
the point $(1,1)$ satisfies both these equations
and corresponds to the case when the quadrilateral $ABCD$
is a parallelogram. Thus, the theorem is completely proved.
\end{proof}

\begin{remark}
It is interesting to obtain a variant of Theorem \ref{cr2.1.thm.1} for $k\in (-1,0)$.
One need to change the setup a little.
For a given $k\in (-1,\infty)$ and a given convex quadrilateral $ABCD$, let us define the points $A_1,B_1,C_1,D_1$  by the vector equalities
$$
\overrightarrow{AA_1}=\tfrac{k}{k+1}\,\overrightarrow{AB}, \quad
\overrightarrow{BB_1}=\tfrac{k}{k+1}\,\overrightarrow{BC}, \quad
\overrightarrow{CC_1}=\tfrac{k}{k+1}\,\overrightarrow{CD}, \quad
\overrightarrow{DD_1}=\tfrac{k}{k+1}\,\overrightarrow{DA}.
$$
By analogy with the original setup, the straight lines
$AB_1$, $BC_1$, $CD_1$, $DA_1$ form a quadrilateral $KLMN$ (it~should be noted that  $KLMN$ coincides with $ABCD$ for $k=0$).
If $S$ and $s$ the areas of the quadrilaterals $ABCD$ and $KLMN$ respectively, then
the problem (generalized Popoviciu’s problem) is the following: for a given $k\in (-1,\infty)$, to determine sharp bounds for the ratio $s/{S}$
on the set of all non-degenerate convex quadrilaterals $ABCD$.
\end{remark}


\vspace{5mm}

\begin{thebibliography}{99}

\bibitem{Nikk}
Nikonorov~Yu.G.
\emph{Some problems of Euclidean geometry,}
Preprint, Rubtsovsk Industrial Institute, Rubtsovsk, 1998, 32~p.

\bibitem{Shkl}
Shklyarskii D.O., Chentsov N.N., Yaglom, I.M.
\emph{Geometrical estimates and problems from combinatorial geometry} (Russian), Izdat. ``Nauka'', Moscow, 1974, 383 p.


\end{thebibliography}
\end{document}